\documentstyle[11pt]{article}
\begin{document}

\title{Sequences related to the Pell generalized equation}
\author{Mario Catalani\\
Department of Economics, University of Torino\\
Via Po 53, 10124 Torino, Italy\\
mario.catalani@unito.it}
\date{}
\maketitle
\begin{abstract}
We consider sequences of the type
$A_n=6A_{n-1}-A_{n-2},
\quad A_0=r,\, A_1=s$ ($r$ and $s$ integers)
and show that all sequences
that solve particular
cases of the Pell generalized equation
are expressible as a constant times one of four particular
sequences of the same type.
\end{abstract}

\bigskip
\bigskip

\noindent
Let $\alpha=3+2\sqrt{2},\,\beta=3-2\sqrt{2}$ be the roots of the polynomial
$x^2-6x+1$. Note that $\alpha+\beta =6,\,\alpha\beta =1,\,\alpha-\beta
=4\sqrt{2}$.
Also let
$\gamma=1+\sqrt{2},\,\delta=1-\sqrt{2}$. Then
$\gamma^2=\alpha,\,\delta^2=\beta,\,\gamma\delta=-1$. We take
$\gamma=\alpha^{1\over 2},\,\delta=-\beta^{1\over 2}$.

\noindent
Consider the sequence $A_n$ defined by
\begin{eqnarray}
\label{eq:sequenzabase}
A_n=6A_{n-1}-A_{n-2},
\quad A_0=r,\, A_1=s,
\end{eqnarray}
where $r$ and $s$ are integers.

\noindent
The object of this contribution is to show that all sequences
of the type given by Equation~\ref{eq:sequenzabase} that solve particular
cases of the Pell generalized equation (see \cite{robertson})
are expressible as a constant times one of four particular
sequences of the same type.

\noindent
The generating function of $A_n$ is given by
$$g(x)={r+(s-7r)x+(6r-s)x^2\over (1-x)(1-6x+x^2)},$$
from which we get the closed form
$$A_n={2s\gamma\alpha-2r\gamma\over 8\sqrt{2}\gamma^3}
\alpha^n -{2r\gamma\beta+2s\delta\beta\over 8\sqrt{2}\delta^3}\beta^n.$$

\noindent
Two important particular cases of Equation~\ref{eq:sequenzabase}
are the sequences
\begin{equation}
\label{eq:uno}
T_n={\alpha^n-\beta^n\over 4\sqrt{2}},
\end{equation}
\begin{equation}
\label{eq:due}
L_n=\alpha^n+\beta^n.
\end{equation}
The initial conditions are, respectively,
$T_0=0,\,T_1=1,$ and
$L_0=2,\,L_1=6$.
They are, respectively, sequences A001109 and A003499 in \cite{sloane}.
$T_n$ is related to triangular numbers: see \cite{ram}.
The relationships between $T_n$ and $L_n$ are of the same genre as those
between Fibonacci and Lucas numbers: so a wealth of known identities
relating Fibonacci and Lucas numbers translates to our pair.

\noindent
Now we establish some identities concerning $\alpha$ and $\beta$ that allow
us to introduce the other sequences necessary for our argument. First
note that
from $\alpha+\beta+2\sqrt{\alpha\beta}=8$ it follows
$\left (\alpha^{1\over 2}+\beta^{1\over 2}\right )^2=8$, that is
$\alpha^{1\over 2}+\beta^{1\over 2}=2\sqrt{2}$. Also
from $\alpha+\beta-2\sqrt{\alpha\beta}=4$ it follows
$\left (\alpha^{1\over 2}-\beta^{1\over 2}\right )^2=4$, that is
$\alpha^{1\over 2}-\beta^{1\over 2}=2$.
Then
\begin{enumerate}
\item
$$\left (\alpha^{n+{1\over 2}}+\beta^{n+{1\over 2}}\right )
\left (\alpha^{1\over 2}-\beta^{1\over 2}\right )=
\alpha^{n+1}-\alpha^n+\beta^n-\beta^{n+1},$$
from which we get
$${\alpha^{n+{1\over 2}}+\beta^{n+{1\over 2}}\over 2\sqrt{2}}=
T_{n+1}-T_n=B_n,$$
which is sequence A001653 in \cite{sloane}. We have $B_0=1,\,B_1=5$.
\item $$\left (\alpha^{n+{1\over 2}}+\beta^{n+{1\over 2}}\right )
\left (\alpha^{1\over 2}+\beta^{1\over 2}\right )=
\alpha^{n+1}+\alpha^n+\beta^n+\beta^{n+1}.$$
Then
\begin{eqnarray*}
\sqrt{2}\left (\alpha^{n+{1\over 2}}+\beta^{n+{1\over 2}}\right )&=&
{\alpha^{n+1}+\alpha^n+\beta^n+\beta^{n+1}\over 2}\\
&=&{L_{n+1}+L_n\over 2}\\
&=&E_n,
\end{eqnarray*}
which is sequence A077445  in \cite{sloane}. We have $E_0=4,\,E_1=20$.
We see that $E_n=4B_n$.
\item
$$\left (\alpha^{n+{1\over 2}}-\beta^{n+{1\over 2}}\right )
\left (\alpha^{1\over 2}+\beta^{1\over 2}\right )=
\alpha^{n+1}+\alpha^n-\beta^n-\beta^{n+1}.$$
From this we get
$$\alpha^{n+{1\over 2}}-\beta^{n+{1\over 2}}=2(T_{n+1}+T_n)=C_n.$$
This is sequence A077444 in \cite{sloane}. We have $C_0=2,\,C_1=14$.
We can write as well
$$C_n= \gamma^{2n+1}+\delta^{2n+1}.$$
Then ${C_n\over 2}$ are the NSW numbers (sequence A002315 in \cite{sloane},
also see \cite{primeglossary}).
Also ${C_n\over 2}=R_{2n+1}$, where the $R_n$ sequence gives the
numerators of continued fraction convergents to $\sqrt{2}$
(sequence A001333 in \cite{sloane}).

\item
$$\left (\alpha^{n+{1\over 2}}-\beta^{n+{1\over 2}}\right )
\left (\alpha^{1\over 2}-\beta^{1\over 2}\right )=
\alpha^{n+1}-\alpha^n-\beta^n+\beta^{n+1}.$$
From this we get
$$\alpha^{n+{1\over 2}}-\beta^{n+{1\over 2}}={L_{n+1}-L_n\over 2}=C_n.$$
\end{enumerate}

\noindent
Now we are going to study in more details the $L_n$ sequence.
In the $L_n$ sequence there are hidden some sequences of squares.
More precisely, we have that $L_{2n}+2,\, 2(L_{2n}-2),\,L_{2n+1}-2,\,
2(L_{2n+1}+2)$ are perfect squares.
\begin{enumerate}
\item
\begin{eqnarray*}
L_n^2&=&(\alpha^n+\beta^n)^2\\
&=&\alpha^{2n}+\beta^{2n}+2\\
&=&L_{2n}+2.
\end{eqnarray*}
\item
\begin{eqnarray*}
(8T_n)^2&=&64\left ({\alpha^n-\beta^n\over 4\sqrt{2}}\right )^2\\
&=&2(\alpha^{2n}+\beta^{2n}-2)\\
&=&2(L_{2n}-2).
\end{eqnarray*}
\item
\begin{eqnarray*}
L_{2n+1}-2&=&\alpha^{2n+1}+\beta^{2n+1}-2\\
&=&\left (\alpha^{n+{1\over 2}}-\beta^{n+{1\over 2}}\right )^2\\
&=&C_n^2.
\end{eqnarray*}
\item
\begin{eqnarray*}
2(L_{2n+1}+2)&=&2(\alpha^{2n+1}+\beta^{2n+1}+2)\\
&=&\left (\sqrt{2}\left (\alpha^{n+{1\over 2}}+\beta^{n+{1\over 2}}\right
)\right )^2\\
&=&E_n^2.
\end{eqnarray*}
\end{enumerate}
Also we obtain easily
\begin{equation}
4(8T_n^2+1)=L_{2n}+2=L_n^2,
\end{equation}
\begin{equation}
2C_n^2+8=2(L_{2n+1}+2)=E_n^2,
\end{equation}
and finally
\begin{equation}
4(2B_n^2-1)=L_{2n+1}-2=C_n^2,
\end{equation}
which means that $8T_n^2+1,\,2C_n^2+8$ and $2B_n^2-1$ are perfect squares.

\noindent
After some tedious algebra we can write the following formula involving
the square of $A_n$:
\begin{equation}
\label{eq:fondamentale}
32A_n^2+2(r^2+s^2-6rs)=(r^2-s^2)L_{2n-2}+(6s^2-2rs)L_{2n-1}.
\end{equation}
We would like to find conditions on $r$ and $s$ under which the LHS can
be reduced to a perfect square. This happens only if we can find values
of $r$ and $s$ for which the RHS can be
written as a constant times $L_n$ for some $n$: this because of the results
obtained before for $L_n$. In this case we would have solutions to
particular generalized Pell equations.

\noindent
The first elementary cases arise when we set equal to zero one of the
two coefficients on the RHS. If we set $r^2=s^2$ we get $s=\pm r$.
If $s=r$ then $A_0=r,\,A_1=r,\,A_2=5r,\,A_3=29r,\,\ldots$ so we can conclude
$A_n=rB_{n-1}$. If $s=-r$
then $A_0=r,\,A_1=-r,\,A_2=-7r,\,A_3=-41r,\,\ldots$ so we can conclude
$A_n=-r{C_{n-1}\over 2}$.

\noindent
If we set $6s^2-2rs=0$ we obtain $s=0$ or $3s=r$. If $s=0$
then $A_0=r,\,A_1=0,\,A_2=-r,\,A_3=-6r,\,\ldots$ so we can conclude
$A_n=-rT_{n-1}$. If $3s=r$, setting $r=3k$ we have
$A_0=3k,\,A_1=k,\,A_2=3k,\,A_3=17k,\,\ldots$ so that
$A_n=k{L_{n-1}\over 2}$.

\noindent
Now let $a_0=2rs-6s^2,\,a_1=r^2-s^2$ and define the
recurrence
$$a_n=6a_{n-1}-a_{n-2},\qquad n\ge 2.$$
Then
\begin{eqnarray*}
32A_n^2+2(r^2+s^2-6rs)&=&a_1L_{2n-2}-a_0L_{2n-1}\\
&=&a_1(6L_{2n-1}-L_{2n})-a_0L_{2n-1}\\
&=&(6a_1-a_0)L_{2n-1}-a_1L_{2n}\\
&=&a_2L_{2n-1}-a_1L_{2n}\\
&=&a_2(6L_{2n}-L_{2n+1})-a_1L_{2n}\\
&=&(6a_2-a_1)L_{2n}-a_2L_{2n+1}\\
&=&a_3L_{2n}-a_2L_{2n+1}\\
&=& \cdots\quad\cdots\quad\cdots
\end{eqnarray*}
So in general we have
$$32A_n^2+2(r^2+s^2-6rs)=a_{m+3}L_{2n+m}-a_{m+2}L_{2n+m+1},\quad
m=-2,\,-1,\,0,\,1 ...$$
Our problem is solved if for some integer $h$
$$a_{m+2}=-6h=-T_2h,\qquad a_{m+3}=-h=-T_1h,$$
since then
\begin{eqnarray*}
32A_n^2+2(r^2+s^2-6rs)&=&-hL_{2n+m}+6hL_{2n+m+1}\\
&=&hL_{2n+m+2}.
\end{eqnarray*}
Now
\begin{eqnarray*}
a_{m+1}&=&6a_{m+2}-a_{m+3}\\
&=&-6T_2h+T_1h\\
&=&-T_3h,
\end{eqnarray*}
\begin{eqnarray*}
a_m&=&6a_{m+1}-a_{m+2}\\
&=&-6T_3h+T_2h\\
&=&-T_4h,
\end{eqnarray*}
and so on. Then the conditions are
$$a_0=-T_{m+4}h,\qquad a_1=-T_{m+3}h.$$
For sake of simplicity, write $T_{m+4}=t_1,\,T_{m+3}=t_0$. So we have
\begin{equation}
\label{eq:uno}
2rs-6s^2+t_1h=0,
\end{equation}
\begin{equation}
\label{eq:due}
r^2-s^2+t_0h=0.
\end{equation}
Solving for $r$ Equation~\ref{eq:uno} and solving for $h$
Equation~\ref{eq:due} after insertion of the value for $r$ we get
$$h={2(-s^2t_0+3s^2t_1\pm s^2)\over t_1^2}.$$
In obtaining this result we used the identity
$$T_n^2-6T_nT_{n+1}+T_{n+1}^2=1.$$
This can be proved in the following way:
$$T_n^2={\alpha^{2n}+\beta^{2n}-2\over 32},\quad
T_{n+1}^2={\alpha^{2n+2}+\beta^{2n+2}-2\over 32},$$
$$T_nT_{n+1}={\alpha^{2n+1}+\beta^{2n+1}-6\over 32}.$$
Then
$$T_n^2-6T_nT_{n+1}+T_{n+1}^2=
{1\over 32}(L_{2n}+L_{2n+2}-6L_{2n+1}+32),$$
but
$$L_{2n}+L_{2n+2}-6L_{2n+1}=0,$$
from which the desired identity.

\noindent
If in $h$ we take the plus sign and insert into the expression for $r$
we get
\begin{equation}
\label{eq:meno}
r={s(t_0-1)\over t_1}.
\end{equation}
On the other hand if we take the minus sign
we obtain
\begin{equation}
\label{eq:piu}
r={s(t_0+1)\over t_1}.
\end{equation}
Next we will prove the following identities
\begin{equation}
\label{eq:ideuno}
T_{2n+1}B_{n-1}=(1+T_{2n})B_n,
\end{equation}
\begin{equation}
\label{eq:idedue}
T_{2n}L_{n-1}=(1+T_{2n-1})L_n.
\end{equation}
Identity~\ref{eq:ideuno}:
\begin{eqnarray*}
T_{2n+1}B_{n-1}&=&{\alpha^{2n+1}-\beta^{2n+1}\over 4\sqrt{2}}
{\alpha^{n-{1\over 2}}+\beta^{n-{1\over 2}}\over 2\sqrt{2}}\\
&=&{\alpha^{3n+{1\over 2}}+\alpha^{2n+1}\beta^{n-{1\over 2}}
-\alpha^{n-{1\over 2}}\beta^{2n+1}-\beta^{3n+{1\over 2}}\over 16}\\
&=&{\alpha^{3n+{1\over 2}}+\alpha^{n+1+{1\over 2}}
-\beta^{n+1+{1\over 2}}-\beta^{3n+{1\over 2}}\over 16},
\end{eqnarray*}
where we used several times the fact that $\alpha\beta=1$.
\begin{eqnarray*}
B_n(1+T_{2n})&=&{\alpha^{n+{1\over 2}}+\beta^{n+{1\over 2}}\over 2\sqrt{2}}
+
{\alpha^{3n+{1\over 2}}-\alpha^{n+{1\over 2}}\beta^{2n}
+\alpha^{2n}\beta^{n+{1\over 2}}-\beta^{3n+{1\over 2}}\over 16}\\
&=&{4\sqrt{2}\alpha^{n+{1\over 2}}+4\sqrt{2}\beta^{n+{1\over 2}}\over 16}+\\
&&\quad\quad
+{\alpha^{3n+{1\over 2}}-\alpha^{n+{1\over 2}}\beta^{2n}
+\alpha^{2n}\beta^{n+{1\over 2}}-\beta^{3n+{1\over 2}}\over 16}\\
&=&{\alpha^{n+1+{1\over 2}}-\alpha^{n+{1\over 2}}\beta
+\alpha\beta^{n+{1\over 2}}-\beta^{n+1+{1\over 2}}\over 16}+\\
&&\quad\quad
+{\alpha^{3n+{1\over 2}}-\alpha^{n+{1\over 2}}\beta^{2n}
+\alpha^{2n}\beta^{n+{1\over 2}}-\beta^{3n+{1\over 2}}\over 16}\\
&=&{\alpha^{3n+{1\over 2}}+\alpha^{n+1+{1\over 2}}
-\beta^{n+1+{1\over 2}}-\beta^{3n+{1\over 2}}\over 16},
\end{eqnarray*}
where we used again $\alpha\beta=1$ and $\alpha-\beta-4\sqrt{2}$.
Hence
$$T_{2n+1}B_{n-1}=(1+T_{2n})B_n={C_{3n}+C_{n+1}\over 16}.$$
Identity~\ref{eq:idedue}:
\begin{eqnarray*}
T_{2n}L_{n-1}&=&{\alpha^{2n}-\beta^{2n}\over 4\sqrt{2}}(\alpha^{n-1}+
\beta^{n-1})\\
&=&{\alpha^{3n-1}+\alpha^{2n}\beta^{n-1}-\alpha^{n-1}\beta^{2n}
-\beta^{3n-1}\over 4\sqrt{2}}\\
&=&
{\alpha^{3n-1}+\alpha^{n+1}-\beta^{n+1}
-\beta^{3n-1}\over 4\sqrt{2}}.
\end{eqnarray*}
\begin{eqnarray*}
(1+T_{2n-1})L_n&=&\alpha^n+\beta^n +
(\alpha^n+\beta^n){\alpha^{2n-1}-\beta^{2n-1}\over 4\sqrt{2}}\\
&=&{\alpha^{n+1}-\alpha^n\beta+\alpha\beta^n-\beta^{n+1}\over 4\sqrt{2}}+\\
&&\quad\quad
+{\alpha^{3n-1}-\alpha^n\beta^{2n-1}+\alpha^{2n-1}\beta^n-\beta^{3n-1}\over
4\sqrt{2}}\\
&=&
{\alpha^{3n-1}+\alpha^{n+1}-\beta^{n+1}
-\beta^{3n-1}\over 4\sqrt{2}}.
\end{eqnarray*}
Hence
$$T_{2n}L_{n-1}=(1+T_{2n-1})L_n=T_{3n-1}+T_{n+1}.$$
Along these lines we can prove these two other identities
\begin{equation}
\label{eq:idetre}
(T_{2n}-1)C_n=T_{2n+1}C_{n-1}={B_{3n}-B_{n+1}\over 2},
\end{equation}
\begin{equation}
\label{eq:idequat}
(T_{2n+1}-1)T_{n+1}=T_nT_{2n+2}={L_{3n+2}-L_{n+2}\over 32}.
\end{equation}
Going back to Equation~\ref{eq:meno} if $m=2k$ we get using
Identity~\ref{eq:idequat}
$${s\over r}={T_{k+2}\over T_{k+1}},$$
from which
$$r=\mu T_{k+1}=A_0,\quad s=\mu T_{k+2}=A_1,$$
so that $A_n=\mu T_{n+k+1}.$
If $m=2k+1$ we get using
Identity~\ref{eq:idetre}
$${s\over r}={C_{k+2}\over C_{k+1}},$$
from which
$$r=\mu C_{k+1}=A_0,\quad s=\mu C_{k+2}=A_1,$$
so that $A_n=\mu C_{n+k+1}.$
On the other hand in the case of
Equation~\ref{eq:piu} if $m=2k$ we get using
Identity~\ref{eq:due}
$${s\over r}={L_{k+2}\over L_{k+1}},$$
from which
$$r=\mu L_{k+1}=A_0,\quad s=\mu L_{k+2}=A_1,$$
so that $A_n=\mu L_{n+k+1}.$
Finally if $m=2k+1$ we get using
Identity~\ref{eq:uno}
$${s\over r}={B_{k+2}\over B_{k+1}},$$
from which
$$r=\mu B_{k+1}=A_0,\quad s=\mu B_{k+2}=A_1,$$
so that $A_n=\mu B_{n+k+1}.$

\noindent
The conclusion is that any sequence we looked for is expressible
as a constant times one of the four sequences $L_n$, $T_n$, $B_n$ and
$C_n$.

\end{document}